\documentclass[11pt]
{amsart}
\def\fulldate{1999-01-07 18:56:32}

\title{Antichains in products of linear orders} 

\author{Martin Goldstern}
\address{Institut f\"ur Algebra\\
Technische Universit\"at Wien\\
Wiedner Hauptstra\ss e 8-10/118.2\\
A-1040 Wien, Austria}

\email{Martin.Goldstern@tuwien.ac.at}
\urladdr{http://info.tuwien.ac.at/goldstern/}
\thanks{The first author is partially supported by the 
Austrian Science Foundation, FWF grant P13325-MAT.}

\author{Saharon Shelah}
\thanks{The second author is supported by the
   German-Israeli Foundation for Scientific Research \& Development
   Grant No. G-294.081.06/93.   Publication number 696.}
\address{Department of Mathematics\\
Hebrew University of Jerusalem\\
Givat Ram\\
91904 Jerusalem, Israel}
\email{shelah@math.huji.ac.il}
\urladdr{http://math.rutgers.edu/\char`\~shelah/}

\usepackage{amsfonts}

\def\xtend#1#2#3{\def\comparetoempty{#2}%
        \ifx\comparetoempty\emptyseq {\rm extend[#1,#3]}\else
	        {\rm extend}(#1,#2{\mapsto#3})\fi}
\def\xtend#1#2#3{\def\comparetoempty{#2}%
        \ifx\comparetoempty\emptyseq #1^{(#3)}\else
                                     #1^{(#2\mapsto#3)}\fi}
\def\xtend#1#2#3{\def\comparetoempty{#2}%
\ifcase #3 \def\opposite{1}\or\def\opposite{0}\else\def\opposite{??}\fi	
        \ifx\comparetoempty\emptyseq #1\conc \bar #3\else
                             #1\conc
               { \scriptstyle\{#2 \mapsto #3 \text{ else }\bar\opposite\}}\fi}

		\def\emptyseq{}
\def\<#1>{_{\langle#1\rangle}}
\def\pp#1>{\langle#1\rangle}

\newcommand{\inc}{{\rm inc}}

\newcommand{\lm}{{\lambda_{m^*}^+}}
\renewcommand{\P}{{\mathbb{P}}}
\newcommand{\eeta}{\vec{\eta}}\renewcommand{\eeta}{\bar{\eta}}
\newcommand{\nnu}{\vec{\nu}}\renewcommand{\nnu}{\bar{\nu}}
\newcommand{\nnnu}{\name{\nnu}}
\newcommand{\rrho}{\vec{\rho}}\renewcommand{\rrho}{\bar{\rho}}
\newcommand{\forces}{\Vdash}

\newcommand{\lc}{%
	\mathrel{<\kern-2.2mm\raise0.25mm\hbox{$\scriptstyle\circ$}}}

\newcommand{\XX}{\trianglelefteq}
\newcommand{\YY}{\trianglerighteq}
\newcommand{\conc}{\!^\frown\!}

\newcount\skewfactor
\def\mathunderaccent#1#2 {\let\theaccent#1\skewfactor#2
\mathpalette\putaccentunder}
\def\putaccentunder#1#2{\oalign{$#1#2$\crcr\hidewidth
\vbox to.2ex{\hbox{$#1\skew\skewfactor\theaccent{}$}\vss}\hidewidth}}
\def\name{\mathunderaccent\widetilde-3 }

\let\oldref\ref
\let\ref\oldref

\usepackage{amssymb}
\usepackage{amsfonts}
\usepackage{amsthm}
\usepackage{amsmath}

\def\xx#1 {\newtheorem{#1}[thm]{#1}}
\xx Theorem
\xx Lemma
\xx {Main Lemma}
\xx Corollary
\xx {Fact A}
\xx {Fact B}
\xx Conclusion
\theoremstyle{definition}
\xx Notation
\xx {Definition and Notation}
\xx Definition
\xx Setup
\xx Remark
\xx Note
\xx {Proof sketch}
\xx Claim
\xx Fact
\xx {Fact and Definition}
\xx {Definition and Fact}
\xx Construction
\xx {Remark and Definition}

\newcommand{\on}{{\restriction}}
\def\pre#1^#2{{}^{#2}#1}

\def\ii{\begin{ENUMERATE}}
\def\jj{\end{ENUMERATE}}

\def\texcases#1{\left \{\,\vcenter {\normalbaselines \m@th \ialign
{$##\hfil $&\quad ##\hfil \crcr #1\crcr }}\right.}

\def\itm#1 {\item[(#1)]}
\def\conc{{}^\frown}

\begin{document}

\begin{abstract}
We show that: \\
\begin{enumerate}
\item For many regular cardinals $ \lambda$ (in particular, for all successors
   of singular strong limit cardinals, and for all successors of
   singular $\omega $-limits), for all $n\in \{2,3,4,\ldots\}$:
   There is a linear order $L$ such that $L^n$ has no
   (incomparability-)antichain of 
   cardinality $\lambda$, while $L^{n+1}$ has an antichain of
   cardinality $\lambda$. 
\item For any nondecreasing sequence $(\lambda_n: n \in \{2,3,4,\ldots
   \})$ of infinite cardinals it is consistent that there  is a linear
   order $L$ such that $L^n$ has an antichain of cardinality
   $\lambda_n$, but not one of cardinality $\lambda_n^+$. 
\end{enumerate}
\end{abstract}

\maketitle

%  \renewcommand{\theenumi}{$\bullet_{\arabic{enumi}}$}
%\renewcommand{\labelenumi}{$\bullet_{\theenumii}$}

%\hrule\bigskip \begin{center}\Large DRAFT \end{center}\bigskip\hrule\bigskip

\centerline{\fulldate}

\pagestyle{myheadings}
\newcommand{\msection}[1]{%\newpage
\section{#1}
\markboth{\fulldate.\hfil  #1\hfil}{\fulldate.\hfil  #1\hfil}}

\markboth{\fulldate}{\fulldate}

\section{Introduction}\label{section0}

\begin{Definition} For any partial ordering $(P,\le)$ define $\inc(P)$ as 
$$ \inc(P) = \sup \{ |A|^+ : A \subseteq P \text{ is an antichain}\}$$
Here, an antichain is a set of pairwise {\bf incomparable} elements. 

In other words, $ \kappa < \inc(P)$ iff there is an antichain of
cardinality $\kappa $. 
\end{Definition}

Haviar and Ploscica in \cite{HP98} asked:  Can there be a linear ordering
$L$ such that $\inc(L^n)\not=\inc(L^k)$ for some natural numbers
$k$.  (Here, $L^n$ and $L^k$ carry the product, or pointwise, order.)

Farley \cite{farley} has pointed out that for any singular cardinal
$\kappa$ there is a linear order $L$ of size $\kappa$ such that 
$\inc(L^2) = \kappa$, $\inc(L^3) = \kappa^+$. 

So we will be mainly interested in this question for {\em regular}
cardinals. 
  First we
show in ZFC that there are many successor cardinals $\lambda $ (including
$\aleph_{\omega+1}$)  with the following
property:  
\begin{quote} For every $n>2$ there is a linear  order $L$ of size
$\lambda $ such that $\inc(L^n)\le\lambda$, $\inc(L^{n+1}) =
\lambda^+$. 
\end{quote}
This proof is given in section \ref{section1}.  It uses a basic fact from pcf
theory.

% BUT COULD BE REWRITTEN TO AVOID THAT??

We then show that there are (consistently) many possible behaviours for
the sequence $(\inc(L^n):n=2,3,4,\ldots)$. More precisely, we show
that for any nondecreasing sequence of  infinite  cardinals $(\lambda_n:
2\le n < \omega)$ there is a cardinal-preserving extension of the
universe in which we can find a linear order $L$ such that for all $n\in
\{2,3,\ldots\}$: $\inc(L^n)=\lambda_n^+$. 

For example,  it is consistent
that there is a linear order 
 $L$ such that $L^2$ has no uncountable antichain, while $L^3$ does. 

   Here we use forcing. 
 The heart of this second  proof is the
well-known ${\Delta}$-system lemma.   

\section{A ZFC proof}\label{section1}

\newcommand{\D}{{D^{bd}_{\mu} }}

Let ${\mu} $ be a regular cardinal.  We will
write $\D$ for the filter of cobounded sets, i.e., the filter dual to
the ideal $\{A \subseteq {\mu}: \exists i<{\mu} \,\, A \subseteq
i\}$.  

$\prod_{i<{\mu}} \lambda_i$ is the set of all functions $f$ with
domain ${\mu}$ satisfying $f(i) <\lambda_i$ for all~$i$.    The
relation $f \sim_\D g \Leftrightarrow\{i:f(i)=g(i)\}\in \D$ is an
equivalence relation.   We call the quotient structure
 $\prod_i
\lambda_i/\D$
 (and we often do not distinguish between a function $f$ and 
its equivalence class).   
 $\prod_i
\lambda_i/\D$
is partially ordered by the relation 
$$f<_\D g\text{  iff }\{i<{\mu}:
f(i)<g(i) \} \in \D$$ 

For any partial order $(P,{\le})$ and any regular cardinal $\lambda $
we say $\lambda = tcf(P)$ (``$\lambda$ is the true cofinality of
$P$'') iff there is an increasing sequence $(p_i: i<\lambda)$ such that
$ \forall p\in P \, \exists i<\lambda : p \le p_i$.

\begin{Theorem}\label{maintheorem} Assume that
\begin{enumerate}
\item  $(\lambda_i:i<{\mu})$ is an increasing sequence of regular
cardinals
\item For each $j< \lambda$, $\bigl|\prod_{i<j} \lambda_i \bigr| < \lambda_j$
\item $\lambda $ is regular and 
$tcf(\prod \lambda_i \big/ \D) = \lambda $, 
\item $n \ge 2$.

\end{enumerate}
Then there is a linear order $J$ of size $\lambda$ such that 
\begin{itemize}
\item $J^{n+1}$ has an antichain of size $\lambda$
\item $J^{n}$ has no antichain of size $\lambda$
\end{itemize}
\end{Theorem}

\begin{Remark} Whenever $\lambda = \mu^+ $ is the successor of a singular
cardinal  $\mu$ such that 
\begin{enumerate}
\item Either $\mu$ is a strong limit cardinal
\item or at least $ \forall \kappa < \mu: \kappa^{<cf(\mu)}<\mu$
\end{enumerate}
 then we can find a sequence $(\lambda_i:i<cf(\mu))$ as above.
For example, if $\lambda = \aleph_{\omega+1}$, then there is an
increasing sequence $(n_k: k \in \omega )$ of natural numbers such
that  $tcf(\prod_{k\in \omega } \aleph_{n_k} /D^{bd}_\omega ) =
\aleph_{\omega+1}$.  See \cite[??]{card}.

\end{Remark}

% REMARK:   UNBOUNDED FILTER IS NO LOSS OF GENERALITY< SINCE OTHER CASES
% ARE REDUCED TO IT??

The proof of theorem~\ref{maintheorem}
 will occupy the rest of this section.  We fix a sequence
$(f_\alpha: \alpha < \lambda)$ as in the assumption of the theorem. 
 We start by writing
${\mu} = \bigcup_{\ell=0}^n A_\ell$ as a disjoint union of~$n+1$
many $\D$-positive (i.e., unbounded) sets. 
 For $\ell=0,\ldots, n$ we define a linear
order $<_\ell$ on~$\lambda$ as follows: 

\begin{Definition}

For any two functions $f,g\in \prod_i \lambda _i $ we define 
\begin{equation}
d(f,g) = \sup \{i<{\mu}: f\on i = g \on i \}
= \max \{i<{\mu}: f\on i = g \on i \}
\label{0}
\end{equation}

That is, if $f \not= g$ we have that $d(f,g) = \min
 \{j: f(j) \not=g (j)\}$ is the first point where
$f$ and $g$ diverge.

For $\alpha, {\beta}\in \lambda$ we define $\alpha< _\ell {\beta} $ iff:
\begin{equation}
\begin{split}\label{1}
\text{letting $i_{\alpha,\beta}:= d(f_\alpha, f_\beta)$},& \\
\text{either } &i_{\alpha,\beta}\in A_\ell \text{  and }
  f_\alpha (i_{\alpha,\beta})<
  f_\beta(i_{\alpha,\beta})  \\
\text{or }   & i_{\alpha,\beta}\notin A_\ell
        \text{ and }
          f_\alpha (i_{\alpha,\beta})> f_{\beta}(i_{\alpha,\beta})
\end{split}
\end{equation}
\end{Definition}

We now define $J$ to be the ``ordinal sum'' of all the orders
$<_\ell$:
\begin{Definition}
Let $$ J = \bigcup_{\ell=0}^n \{\ell\} \times (\lambda, <_\ell)$$
with the ``lexicographic'' order, i.e., $ \pp\ell_1,\alpha_1 > <  \pp\ell_2,
\alpha_2 >$ iff $\ell_1<\ell_2$, or $\ell_1=\ell_2$ and $\alpha_1
<_{\ell_1} \alpha_2$. 
\end{Definition}

\begin{Claim} $J^{n+1} $ has an antichain of size $\lambda$. 
\end{Claim}
\begin{proof} Let $\vec t_\alpha = ( (0,\alpha), \ldots, (n,\alpha))\in
J^{n+1}$. 

For any $\alpha \not= \beta $ we have to check that $\vec t_\alpha$ and
 $\vec t_\beta $ are incomparable.   Let $i^* = d(f_\alpha, f_\beta)$,
 and find $\ell^*$ such that $i^*\in A_{\ell^*}$.  Wlog assume
 $f_\alpha(i^*)<f_\beta(i^*)$. Then $\alpha <_{\ell^*} \beta$, but
 $\alpha >_{\ell} \beta  $ for all $\ell\not= \ell^*$, i.e., 
$ \pp\ell^*, \alpha > <_J  \pp\ell^*, \beta  >$, but 
$ \pp\ell, \alpha > >_J  \pp\ell, \beta  >$ for all $\ell\not= \ell^*$. 

\end{proof}

\begin{proof}[Proof of \ref{maintheorem}]
It remains  to show that $J^n$ does not have an antichain of size
$\lambda$.  Towards a contradiction, assume that 
 $(\vec t_\beta: \beta < \lambda)$ is an antichain in~$J^m$, $m\le n$,
and $m$ as small as possible.   Let
$\vec t_\beta = (t_\beta(1), \ldots, t_\beta(m))\in J^m$. 
For $k=1,\ldots, m$ we
 can find functions  $\ell_k$, $\xi_k$ such that  
\begin{equation*}
 \forall \beta < \lambda \ \forall k : \ t_\beta(k) =
 \pp\ell_k(\beta), \xi_k(\beta) >
% \label{3}
\end{equation*}

Thinning out we may assume that the functions 
  $\ell_1,\ldots , \ell_n$ are constant. We will again write $\ell_1,
  \ldots, \ell_n$ for those  constant values.

 We may also
assume that for each $k$ the function $ \beta \mapsto \xi_k(\beta)$ is
either constant or strictly increasing.       If any of the
functions   $ \xi_k$  is constant we get a contradiction to the minimality of
$m$, so  all
the $\xi_k$ are strictly increasing.   So we may moreover assume that
$\beta < \gamma $ implies $\xi_k(\beta) < \xi_{k'}(\gamma)$ for all
$k,k'$, and in particular $ \beta \le \xi_k(\beta)$ for all $\beta, k$. 

Now define $g^+_\beta, g^-_\beta \in \prod_{i<\mu} \lambda_i$ for
every $\beta < \lambda$ as follows: 
\begin{equation}\label{4}
\begin{split}
 g_\beta^+(i) &= \max(f_{\xi_{1}(\beta)}(i), \ldots,
f_{\xi_{n}(\beta)}(i)) \\
g_\beta^-(i) &= \min(f_{\xi_{1}(\beta)}(i), \ldots,
f_{\xi_{n}(\beta)}(i)) 
\end{split}
\end{equation}

\subsubsection*{Subclaim}
The set 
\begin{equation}\label{5}
C:= \{i< \mu :\  \forall \beta \,\, \{g_\gamma^-(i): \gamma>\beta\} 
	\text{ is unbounded in~$\lambda_i$}\}
\end{equation}
is in the filter~$\D$, i.e., its complement
$$ S:= \{i< \mu : \exists \beta<\lambda \, \exists s<\lambda_i\,\,
\  \{g^-_\gamma(i): \gamma> \beta\} \subseteq s \}$$
is in  the ideal dual to  $\D$ (i.e.,  is a bounded set). 

\begin{proof}[Proof of the subclaim]
  For each $i\in S$ let $\beta_i < \lambda$ and $h(i)
 < \lambda_i$ be  such that    $\{g_\gamma^-(i): \gamma > \beta_i \}
 \subseteq h(i)$.   Let $\beta^* = \sup\{\beta_i: i \in S\}<\lambda $, and
 extend $h$ arbitrarily to a total function on~$\mu$.  Since the
 sequence $(f_\alpha: \alpha < \lambda)$ is cofinal in~$\prod_i
 \lambda_i/\D$, we can find  $\gamma >
 \beta^*$ such that  $h <_\D f_\gamma$. 

We have $\gamma \le \xi_{k}(\gamma)$ for all~$k$, so the sets 
$$X_k:= \{ i < \mu : h(i) < f_{\xi_{k}(\gamma)}(i)\}$$ are all in
$\D$.  Now if $S$ were positive mod~$\D$, then we could find $i_*\in S
\cap X_1 \cap \cdots \cap X_n$. But then  $i^* \in X_1 \cap \cdots \cap
X_n$ implies $$ h(i_*) < g_\gamma^-(i_*),$$
 and $i\in S$  implies 
$$ g_\gamma^-(i_*) < h(i_*),$$
a contradiction. 

This shows that $C$ is indeed a set in the filter~$\D$. 
\end{proof}

We will now use the fact that  $m < n+1$.   Let
$$\ell^* \in \{0,\ldots, n\} \setminus \{\ell_1,\ldots, \ell_n\}.$$
Since  $A_{\ell^*}$ is positive mod~$\D$,  we can pick 
\begin{equation}\label{7}
i^* \in A_{\ell^*} \cap C
\end{equation}

Using the fact that $i^*\in C$ and definition~(\ref{5}) we 
can find a sequence $(\beta_\sigma: \sigma < \lambda_{i^*})$ such
 that  
\begin{equation}\label{8}
 \forall \sigma < \sigma ' < \lambda_{i^*}: \ 
g^+_{\beta_\sigma}(i^*) < g^-_{\beta_{\sigma'}}(i^*), 
\end{equation}
   We now restrict our attention from 
$(\vec t_{\beta}: \beta  < \lambda)$ to the subsequence 
$(\vec t_{\beta_\sigma}: \sigma < \lambda_{i^*})$; we will
show that this sequence cannot be an antichain.   For notational
simplicity only we will assume $\beta_\sigma = \sigma$ for all $\sigma
< \lambda_{i^*} $. 

 Recall that $\vec t_\sigma = (   \pp\ell_1, \xi_1(\sigma) >, \ldots, 
 \pp\ell_n, \xi_n(\sigma) >)$. 
For each $\sigma < \lambda_{i^*} $ define 
$\vec x_\sigma:= (f_{\xi_{1}(\sigma)}\on i^*,  \ldots, 
f_{\xi_{n}(\sigma)}\on i^*)
\in \prod_{j<i^*} \lambda_j$.
  Since $\bigl|\prod_{j < i^* } \lambda_j\bigr|
< \lambda_{i^*}$, there are only $<\lambda_{i^*} $ many possible values for
$\vec x_\sigma$, so we can find ${\sigma_1}<\sigma_2<\lambda_{i^*} $
 such that  $\vec x_{\sigma_1} = \vec x_{\sigma_2}$. 

Now note that by (\ref{4}) and (\ref{8}) we have 
\begin{equation}\label{9}
 f_{\xi_{k}({\sigma_1})}({i^*} )  \le g^+_{\sigma_1}({i^*}) < 
g^-_{\sigma_2}({i^*}) \le 
 f_{\xi_{k}(\sigma_2)}({i^*} ). 
\end{equation}

Hence $d(f_{\xi_{k}({\sigma_1})}, 
 f_{\xi_{k}(\sigma_2)}) = i^*$  for~$k=1,\ldots, n$. 

Since $i^*\in A_{\ell^*}$ we have for all $k$: $i^* \notin A_{\ell_k}
$ . {}From (\ref{0}), (\ref{1}), (\ref{9}) we get 
$$ \xi_{k}({\sigma_1}) <_{\ell_k} \xi_{ k}(\sigma_2) 
\qquad \text{ for~$k=1,\ldots, n$. }$$
Hence $ \pp\ell_k, \xi_{k}({\sigma_1})  > <  \pp\ell_k, \xi_{
k}(\sigma_2)  >  $ for all
$k$, which means $\vec t_{\sigma_1} < \vec t_{\sigma_2}$. 
\end{proof}

\newpage
\section{Consistency}\label{section2}

\begin{Theorem}\label{3.1}
Assume $ {\aleph_0} \le  \lambda _2 \le \lambda_3 \le \cdots $,
$\lambda_n^{<\kappa } \le \lambda_n$.    
$\kappa ^{<\kappa} = \kappa $. 
Then there is a forcing notion $\P$ 
which satisfies the  $\kappa$-cc and is  $\kappa$-complete, 
and a $\P$-name $\name I$ such that  
$$ \forces_P \name I \subseteq  2^\kappa,  \alpha(\name I^n) = \lambda_n^+$$
\end{Theorem}

\begin{Remark}  At first reading, the reader may want to consider the
special case $\kappa = \aleph_0$, $\lambda_{n+2}= \aleph_n$. 
\end{Remark}

\begin{Notation}\label{xtend}
\begin{enumerate}
\item We let $\lambda_1 =0$, $\lambda_\omega = \sup\{\lambda_n: n <
\omega\}$.
\item It is understood that  $2^\alpha $ is linearly ordered
lexicographically,  and $(2^\alpha)^m$ is partially ordered by the
pointwise order. 
\item 
 For $\alpha \le \beta \le \kappa$, $\eta\in
	2^\alpha$, $\nu \in 2^\beta$, we define
$$ \eta \XX \nu \text{ \  \ iff \ \  $\nu$ extends~$\eta$, i.e., $\eta
	\subseteq \nu$}$$
\item
For $\eeta\in (2^\alpha)^n$, $\eeta = (\eta(0), \ldots, \eta(n-1))$, 
 $\nnu \in  (2^\beta )^n$,  $\nu = (\nu(0), \ldots, \nu(n-1))$, 
 we let
$$\eeta\XX \nnu\text{ \  iff \ }\eta(0) \XX \nu(0), \ldots, \eta(n-1)\XX
	\nu(n-1).$$
\item For $\eta\in 2^\alpha $, $i\in \{0,1\}$ we write $\eta\conc i$
for the element $\nu\in 2^{\alpha + 1}$ satisfying $\eta \XX \nu$,
$\nu(\alpha)=i$. 
\end{enumerate}
\end{Notation}

\begin{Definition}
Let $\eeta\in (2^\alpha)^m$, $k \in \{0,\ldots, {m-1}\}$, $m \ge 2$.
We define 
 $\xtend\eeta{}{1}$,
 $\xtend\eeta{}{0}$,
 $\xtend\eeta{k}{1}$, 
 $\xtend\eeta{k}{0}$
 in~$(2^{\alpha+1})^m$ as follows: All four 
are $\XX$-extensions of~$\eeta$, and: 
\begin{itemize}
\item[--] $\xtend\eeta{}0(n) = \eta(n)\conc 0$ for all~$n<m$.
\item[--] $\xtend\eeta{}1(n) = \eta(n)\conc 1$ for all~$n<m$.
\item[--] $\xtend\eeta k 0(n) = \eta(n)\conc 1$ for all $n\not=k$,
 	  $\xtend\eeta k 0(k) = \eta(n)\conc 0$.
\item[--] $\xtend\eeta k 1(n) = \eta(n)\conc 0$ for all $n\not=k$,
 	  $\xtend\eeta k 1(k) = \eta(n)\conc 1$.
\end{itemize}
\end{Definition}

\bigskip

\setlength{\unitlength}{0.00083300in}%
\begingroup\makeatletter\ifx\SetFigFont\undefined
% extract first six characters in \fmtname
\def\x#1#2#3#4#5#6#7\relax{\def\x{#1#2#3#4#5#6}}%
\expandafter\x\fmtname xxxxxx\relax \def\y{splain}%
\ifx\x\y   % LaTeX or SliTeX?
\gdef\SetFigFont#1#2#3{%
  \ifnum #1<17\tiny\else \ifnum #1<20\small\else
  \ifnum #1<24\normalsize\else \ifnum #1<29\large\else
  \ifnum #1<34\Large\else \ifnum #1<41\LARGE\else
     \huge\fi\fi\fi\fi\fi\fi
  \csname #3\endcsname}%
\else
\gdef\SetFigFont#1#2#3{\begingroup
  \count@#1\relax \ifnum 25<\count@\count@25\fi
  \def\x{\endgroup\@setsize\SetFigFont{#2pt}}%
  \expandafter\x
    \csname \romannumeral\the\count@ pt\expandafter\endcsname
    \csname @\romannumeral\the\count@ pt\endcsname
  \csname #3\endcsname}%
\fi
\fi\endgroup
% \begin{figure}
\begin{picture}(3012,1935)(1201,-2038)
% \thicklines
\put(1276,-286){\line( 0,-1){1200}}
\put(1426,-286){\line( 0,-1){1200}}
\put(1576,-286){\line( 0,-1){1200}}
\put(1726,-286){\line( 0,-1){1200}}
\put(2476,-286){\line( 0,-1){1200}}
\put(2626,-286){\line( 0,-1){1200}}
\put(2776,-286){\line( 0,-1){1200}}
\put(2926,-286){\line( 0,-1){1200}}
\put(3751,-286){\line( 0,-1){1200}}
\put(3901,-286){\line( 0,-1){1200}}
\put(4051,-286){\line( 0,-1){1200}}
\put(4201,-286){\line( 0,-1){1200}}
\put(2736,-211){\makebox(0,0)[lb]{\smash{$0$}}}
\put(2586,-211){\makebox(0,0)[lb]{\smash{$0$}}}
\put(2886,-211){\makebox(0,0)[lb]{\smash{$0$}}}
\put(2436,-211){\makebox(0,0)[lb]{\smash{$0$}}}
\put(4011,-211){\makebox(0,0)[lb]{\smash{$1$}}}
\put(3861,-211){\makebox(0,0)[lb]{\smash{$0$}}}
\put(4161,-211){\makebox(0,0)[lb]{\smash{$0$}}}
\put(3711,-211){\makebox(0,0)[lb]{\smash{$0$}}}
\put(3996,-1686){\makebox(0,0)[lb]{\smash{$k$}}}
\put(2466,-2011){\makebox(0,0)[lb]{$\xtend{\eeta}{}{0}$}}
\put(3711,-2011){\makebox(0,0)[lb]{$\xtend{\eeta}{k}{1}$}}
\put(1431,-2011){\makebox(0,0)[lb]{$\eeta$}}
\end{picture}

% \end{figure}

\begin{Fact}
\begin{enumerate}
\item
  If  $\alpha \le \beta \le \kappa$, 
$\eeta, \eeta'\in (2^\alpha )^n$ are incomparable, $\nnu, \nnu'\in
(2^\beta)^n$,  $\eeta \XX \nnu$, $\eeta' \XX \nnu'$, then also $\nnu $
and $\nnu'$ are incomparable. 
\item
 $\xtend \eeta{}{0} < \xtend \eeta{}{1}$.
\item 
 $\xtend \eeta k 0 $ and $  \xtend \eeta k 1 $ are incomparable. 
\end{enumerate}
\end{Fact}

\begin{Definition}\label{pdef}
We let $\P$ be
the set of all conditions  
$$ p = (u^p,\alpha^p, (\eeta_\xi^p: \xi\in u^p))$$
satisfying the following conditions for all $m$: 
\begin{itemize}
\item[--] $u^p \in [\lambda_\omega]^{<\kappa}$
\item[--] $\alpha^p < \kappa $
\item[--] For all $\xi \in u^p\cap (\lambda_m\setminus \lambda_{m-1})$: 
	$\eeta^p_\xi = (\eta^p_\xi(0),\ldots, \eta^p_\xi(m-1)) \in
	(2^{\alpha^p})^m $.
\item[--] For all $\xi\not= \xi'$ in~$u^p \cap (\lambda_m\setminus
	\lambda_{m-1}) $, $\eeta^p_\xi$ and $\eeta^p_{\xi'}$ are
	incomparable in $(2^{\alpha^p})^m$.
\end{itemize}
	
We define $p \le q $ (``$q$ is stronger than $p$'') iff 
\begin{itemize}
\item[--] $u^p \subseteq  u^q$
\item [--] $\alpha^p \le \alpha^q$
\item[--] for all $\xi\in u^p$, $\eeta^p_\xi\XX \eeta^q_\xi$
\end{itemize}

\end{Definition}

\begin{Fact} \label{3.5}
\begin{enumerate}
\item For all $\alpha < \kappa$: The set $\{p\in \P: \alpha^p\ge \alpha\}$
        is dense in~$\P$.  
\item For all $\xi < \lambda_\omega$: The set $\{p\in \P: \xi \in
        u^p\}$ is dense in~$\P$. 
\end{enumerate}
\end{Fact}

\begin{Fact and Definition}\label{idef}
We let $(\nnnu_\xi: \xi < \lambda_\omega)$ be the ``generic object'',
i.e., a name satisfying
\begin{align*}
\forall  m \in \omega \,\,
\forall p \in \P\,\, 
	\forall \xi \in u^p \cap (\lambda_m \setminus \lambda_{m-1}): \ &
p \forces_\P \nnnu_\xi \in (2^\kappa)^m \\
 \forall p \in \P\,\, 
	\forall \xi \in u^p: \  &
p \forces \eeta^p_\xi \XX \nnnu_\xi
\end{align*}
(This definition makes sense, by fact \ref{3.5}.)

Clearly, $ \forces \xi,\xi'\in \lambda_m \setminus \lambda_{m-1}
        \Rightarrow \nnnu_\xi, \nnnu_{\xi'} 
		\text{  incompatible}$. 

We let  $\forces \name I =
 \bigcup_{m=2}^\infty \{\nu_\xi(\ell): 
		\xi\in \lambda_m\setminus \lambda_{m-1} ,
		\ell< m \}$. 
\end{Fact and Definition}

\begin{Theorem}
Let~$\P$, $\name I$ be as in \ref{pdef} and \ref{idef}.   \\
Then
 $ \forces_\P \inc(\name I^m) = \lambda_m$. 
\end{Theorem}

It is clear that $\P$ is $\kappa$-complete, and $\kappa^+$-cc is proved
by an argument similar to the $\Delta$-system argument below. 
So all  the $\lambda_m $ stay cardinals.

We can show by induction that $ \forces \alpha(\name I^m) > \lambda_m$,
i.e.,  $\name I^m$ has an antichain of size $\lambda_m$:
This is clear if $\lambda_m=\lambda_{m-1}$ (and void if $m=0$);  if 
$\lambda_m>\lambda_{m-1} $ then $(\nnu_\xi: \xi \in
  \lambda_m\setminus \lambda_{m-1})  $ will be forced to be antichain. 

\bigskip

It remains to show that  (for any $m$)
there is no antichain of size $\lambda_m^+$ in
$\name I^m$. 

Fix $m^*\in \omega$, and assume wlog that
$\lambda_{m^*+1}>\lambda_{m^*}$.

[Why is this no loss of generality?   
 If $\lambda_{m^*} = \lambda_\omega$, then the cardinality of~$\name I$ is at most
$\lambda_{m^*}$, and there is nothing to prove.  If $\lambda_{m^*}=
 \lambda_{m^*+1}  < \lambda_ \omega $, then replace $m^*$ by
 $\min\{m \ge m^*: \lambda_m<\lambda_{m+1} \}$]

Towards a contradiction, assume that there is a condition $p $ and a
sequence of names
 $\langle \name \rrho_\beta : \beta <\lm\rangle$ 
such that
  $$ 
p \forces 
\langle \name \rrho_\beta : \beta <\lm\rangle
\text { is an antichain in~$\name I^{m^*}$} $$
Let $\name \rrho_\beta  = (\name \rho_\beta(n): n < m^*)$. For each
$\beta < \lm$ and each $n<m^*$
 we can find a condition $p_\beta \ge p$ and 
$$ {m}(\beta,n)\in \omega 
\qquad 
 \ell(\beta,n)< {m}(\beta,n)
\qquad
\xi_{n}(\beta)\in \lambda_{{m}(\beta,n)} \setminus\lambda_{{m}(\beta,n)-1}
$$
such that  
$$ p_\beta \forces \rho_\beta(n) = \nu_{\xi_{n}(\beta)}(\ell(\beta,n))$$

We will now employ a $\Delta $-system argument.

We define a family $(\zeta^\beta: \beta < \lambda_{m^*}^+)$ of
functions as follows: Let $i_\beta$ be the order type of 
$u^{p_\beta}$, and let 
$$u^\beta = u^{p_\beta} = \{ \zeta^\beta(i): i < i_\beta\}
\ \ \text{in increasing enumeration}$$

By \ref{3.5}.2
 may assume $\xi_{n}(\beta)\in u^\beta$, say $\xi_{n}(\beta) =
\zeta^\beta(i(\beta, n))$. 

By thinning out our alleged antichain 
 $\langle \name \rrho_\beta : \beta <\lm\rangle$ 
 we may assume
\begin{itemize}
\item For some $i^* < \kappa$, for all $\beta$: $i_\beta= i^*$
\item For some $\alpha^*< \kappa$, for all $\beta$:
	$\alpha^{p_\beta}=\alpha^*$
\item For each $i<i^*$ there is some ${m} \<i> $ such that  
  for all $\beta$: $\zeta^\beta(i) \in \lambda_{{m}\<i>} \setminus 
\lambda_{{m}\<i>-1} $
\item For each $i<i^*$ there is some $\eeta\<i>\in
	(2^{\alpha^*})^{{m}\<i>}$ such that  for all $\beta$:
	$\eeta^{p_\beta}_{\zeta^\beta(i)} = \eeta\<i>$. (Here we use
	$\lambda_m^{<\kappa} = \lambda_m$.)
\item the family $\langle u^\beta:\beta< \lm\rangle$ is a
	$\Delta$-system, i.e., 
	there is some set $u^*$ such that for all $\beta \not= \gamma
	$: $u ^\beta \cap u ^ \gamma = u^*$. 
\item Moreover: there is a set $\Delta \subseteq i^*$ such that for
	all $\beta$: $u^* = \{ \zeta^\beta(i): i \in \Delta\}$.
  Since $\zeta^\beta $ is increasing, this also implies
  $\zeta^\beta(i)=\zeta^\gamma(i)$ for~$i\in \Delta$.)

\item The functions
 $(\beta,n) \mapsto \ell(\beta,n)$,
 and $(\beta,n) \mapsto i(\beta,n)$ do not depend on~$\beta$,
 i.e. there are $(\ell_n: n < m^*)$ and $(i_n: n < m^*)$ such that
 $i(\beta,n) = i_n$, $\ell(\beta,n) = \ell_n$ for all~$\beta$. 
\end{itemize}

Note that for~$i\in i^* \setminus \Delta$ all the $\zeta^\beta(i)$ are
	distinct elements of~$\lambda_{{m}\<i>}$, hence:
\begin{quote}
  $i\notin \Delta $ implies $\lambda^+_{m^*} \le \lambda_{{m}\<i>}$, 
	hence  ${m}\<i> > m^*$.
\end{quote}

Now pick $k^* \le  m^*$ such that $k^*\notin \{ \ell_n: n < m^*\}$. Pick
any distinct  $\beta,\gamma< \lambda _{m^*}^+$.
  We will find a condition $q$ extending $p_\beta$
and $p_\gamma$, such that $q \forces \rrho_\beta \le \rrho_\gamma$. 

\newcommand{\dcup}{\mathbin{\dot{\cup}}}

We define $q$ as follows: 
\begin{itemize}
\item   $u^q: = u_\beta \cup u_\gamma = u^* 
\dcup \{\zeta^\beta(i): i \in i^* \setminus \Delta \}
\dcup \{\zeta^\gamma(i): i \in i^* \setminus \Delta \}$.
\item $\alpha^q = \alpha^* + 1 $.
\item For $\xi \in u^*$, say $\xi = \zeta^\beta(i) = \zeta^\gamma(i)$,
recall that $\eeta^{p_\beta}_\xi =  \eeta\<i> = \eeta^{p_\gamma }_\xi
$.   
We let $\eta^q_\xi = \xtend{\eta\<i>}{}{0}$ (see \ref{xtend}). 
\item  For $ \xi =  \zeta^\beta(i)$, $i\in i^* \setminus \Delta$, we
have  $\eeta^{p_\beta}_\xi = \eeta\<i> \in (2^{\alpha^*})^{{m}\<i>}$,
where ${m}\<i> > m^*$.   Hence $ \xtend{\eeta\<i>}{\ell}{1}$ is well-defined. 
We let 
$$ \eeta^q_\xi = \xtend{\eeta\<i>}{\ell^*}{1}$$

\item  For $ \xi =  \zeta^\gamma (i)$, $i\in i^* \setminus \Delta$, we
let 
$$ \eeta^q_\xi = \xtend{\eeta\<i>}{k^*}{0}$$
\end{itemize}

We claim that $q$ is a condition.   The only nontrivial requirement is
the incompatibility of all $\eeta^q_\xi$:  Let $\xi,\xi'\in u^q$,
$\xi\not= \xi'$. 

If $\xi,\xi' \in u^\beta$, then the incompatibility of~$\eeta^q_\xi$
and $\eeta^q_{\xi'}$ follows from the incompatibility of 
$\eeta^{p_\beta}_\xi$
and $\eeta^{p_\beta}_{\xi'}$. The same argument works for~$\xi,\xi'\in
u_\gamma$. 

So let $\xi\in u_\beta\setminus u^*$, $\xi'\in u_\gamma\setminus
u^*$.   Say $\xi = \zeta^\beta(i)$,  $\xi' = \zeta^\gamma (i')$. 
 
If $i\not= i'$, then $\eeta\<i> = \eeta^{p_\beta}_{\zeta ^\beta(i)} 
= \eeta^{p_\gamma }_{\zeta ^\gamma (i)}$  and 	
$\eeta\<i'> =  \eeta^{p_\gamma }_{\zeta ^\gamma (i')}$ are
incompatible. From 
  $\eeta\<i> \XX\eeta^q_\xi$ and 
 $\eeta\<i'> \XX\eeta^q_{\xi'}$ we conclude that also 
$\eeta^q_\xi$ 
and $\eeta^q_{\xi'}$ are incompatible. 

Finally, we consider the case~$i=i'$.

    We have 
$$ \eeta^q_\xi = \xtend{\eeta\<i>}{k^*}{0} \qquad \qquad 
 \eeta^q_{\xi'} = \xtend{\eeta\<i> }{k^*}{1}$$
so by \ref{3.5}.3,
$\eeta^q_\xi$ 
and $\eeta^q_{\xi'}$ are incompatible.

This concludes the construction of~$q$.  We now check that $q \forces
\rrho_\beta \le \rrho_\gamma$, i.e.,
$ q \forces \rho_\beta(n) \le \rho_\gamma(n)$ for all~$n$.  
Clearly, $q \forces \rho_\beta(n) = \nu_{\zeta^\beta(i_n)}(\ell_n) 
\YY \eeta^q_{\zeta^\beta(i_n)} = 
\xtend{\eta\<i_n>}{}0$.   Here we use the fact that $k^* \not=
\ell_n$. Similarly, $q \forces \rho_\gamma (n) = \nu_{\zeta^\beta(i_n)}(\ell_n) 
\YY \xtend{\eta\<i_n>}{}{1}$.    

Hence  $q \forces
\rrho_\beta \le \rrho_\gamma$.    

This concludes the proof of theorem~\ref{3.1}

%\bibliographystyle{lit-unsrt}
% \bibliographystyle{plain}

% \bibliography{listb,other,goldstrn,listx}

\end{document}